\documentstyle[11pt,leqno]{article}
\input amssym.def
\title{On the Spinor Representation of Surfaces in Euclidean 3-Space.
\footnote{Supported by the SFB 288 of the DFG.}}
\date{\today}
\author{Thomas Friedrich, Berlin}
\begin{document}
\maketitle
\section{Introduction}
The Weierstra{\ss} formula describes a conformal minimal immersion of a
Riemann surface $M^2$ into the 3-dimensional Euclidean  space ${\Bbb R}^3$.
It expresses the immersion in terms of a holomorphic function $g$ and a
holomorphic 1-form $\mu$ as the integral

\[ f= \mbox{Re} \left(€\int (1-g^2 , i (1+g^2), 2g) \mu \right) : 
M^2 \to {\Bbb R}^3 . \]

On the other hand let us consider the spinor bundle $S$ over $M^2$. This 
2-dimensional vector bundle splits into

\[ S=S^+ \oplus S^- = \Lambda^0 \oplus \Lambda^{1,0} . \]

Therefore the pair $(g, \mu)$ can be considered as a spinor field $\varphi$ 
on the Riemann surface. The Cauchy-Riemann equation for  $g$ and $\mu$ is 
equivalent to the homogeneous Dirac equation 

\[ D( \varphi) = 0 . \]

The choice of the Riemannian metric in the fixed conformal class of $M^2$ 
is not essential since the kernel of the Dirac operator is a conformal 
invariant.\\

A similar description for {\it an arbitrary} surface $M^2 \hookrightarrow 
{\Bbb R}^3$ is possible and has been pointed out probably  for the first time
by Eisenhardt (1909).  This representation of any surface in ${\Bbb R}^3$ by 
a spinor field $\varphi$ on $M^2$ satisfying the inhomogeneous Dirac equation\\

\mbox{} \hspace{6.5cm} $D(\varphi)= H \varphi$ \hfill $(*)$\\

involving  the mean curvature $H$ of the surface has been used again in some 
recent papers (see [KS], [R], [Tai1], [Tai2], [Tai3]). However,  the mentioned
authors describe the relationship between surfaces in ${\Bbb R}^3$ and 
solutions of the equation $(*)$ in local terms in order to get explicit
formulas. The aim of the present paper is to clarify the mentioned
representation of surfaces in ${\Bbb R}^3$ by solutions of the equation 
$D(\varphi)=H \varphi$ in a {\it geometrically invariant way}.  It turns out 
that the main idea leading to the description of a surface by a spinor 
field $\varphi$ is simple: Consider an immersion $M^2 \hookrightarrow {\Bbb R}^3$ and fix a parallel spinor $\Phi$ on ${\Bbb R}^3$. Then the restriction $\varphi = \Phi_{|M^2}$ of $\Phi$ to the surface is (with respect to the inner geometry of $M^2$) a non-trivial spinor field on $M^2$ and defines a spinor $\varphi^*$ of constant length which is a solution of the inhomogeneous Dirac equation

\[ D(\varphi^*) =H \varphi^* . \]

Conversely, given a solution $\varphi$ of the equation $(*)$ with constant length there exists a symmetric endomorphism $E:T(M^2) \to T(M^2)$ such that the spinor field satisfies a ''twistor type equation''

\[ \nabla_X \varphi =E(X) \cdot \varphi . \]

The resulting integrability conditions for the endomorphism $E$ are exactly the Gau{\ss} and Codazzi equations. As a consequence, the solution $\varphi$ of the Dirac equation

\[ D(\varphi)= H \varphi \quad , \quad |\varphi | \equiv const \, > 0 \]

yields an isometric immersion of $M^2$ into ${\Bbb R}^3$. In a similar way one obtains the description of conformal immersions using the well-known formula for the transformation of the Dirac operator under a conformal change of the metric (see [BFGK]).

\section{The Dirac Operator of a Surface immersed into a Riemannian 3-Manifold.}

Let $Y^3$ be a 3-dimensional oriented Riemannian manifold with a fixed spin structure and denote by $M^2$ an oriented surface isometrically immersed into $Y^3$. Because the normal bundle of $M^2$ is trivial, the spin structure of $Y^3$ induces a spin structure on the Riemannian surface $M^2$. The spinor bundle $S$ of the 3-manifold $Y^3$ yields by restriction the spinor bundle of the surface $M^2$. Over $M^2$ this bundle decomposes into

\[ S=S^+ \oplus S^- \]

where the subbundles $S^{\pm}$ are defined by (see [F])

\[ S^{\pm} = \{ \varphi \in S:  i \cdot e_1 \cdot e_2 \cdot \varphi = \pm \varphi \} \]

Here $\{e_1, e_2 \}$ denotes an oriented orthonormal frame in $T(M^2)$ and $X \cdot \varphi$ means the Clifford multiplication of a spinor $\varphi \in S$ by a vector  $X \in T(M^2)$. Since in the 3-dimensional Clifford algebra the relation 

\[ e_1 \cdot e_2 = e_3 \]

holds, we can replace the Clifford product $e_1 \cdot e_2$ by the normal vector $\vec{N}$ of $M^2 \hookrightarrow Y^3$:

\[ S^{\pm} = \{ \varphi \in S : i \cdot \vec{N} \cdot \varphi = \pm \varphi \} . \]

Consider a spinor field $\Phi$ defined on the 3-manifold $Y^3$. Its restriction $\varphi = \Phi_{|M^2}$ is a spinor field defined on $M^2$ and decomposes  therefore  into $\varphi = \varphi^+ + \varphi^-$ with

\[ \varphi^+ = \frac{1}{2} ( \varphi + i \vec{N} \cdot \varphi ) \quad , \quad \varphi^- = \frac{1}{2} ( \varphi - i \vec{N} \cdot \varphi ) . \]

We denote by $\nabla^{Y^3}$ and $\nabla^{M^2}$ the covariant derivatives in the spinor bundles on $Y^3$ and $M^2$ respectively. For any vector $X \in T(M^2)$ we have  the well-known formula (see [BFGK]) 

\[ \nabla^{Y^3}_X (\Phi) = \nabla^{M^2}_X (\varphi) - \frac{1}{2} (\nabla_X \vec{N}) \cdot \vec{N} \cdot \varphi  . \]

The vector $(\nabla_X \vec{N})$ coincides with the second fundamental form  II: $T(M^2) \to T(M^2)$ of the submanifold $M^2 \hookrightarrow Y^3$. Since II is symmetric the Clifford product $e_1 \cdot \mbox{II}( e_1) + e_2  \cdot \mbox{II}(e_2)$ is a scalar and equals $(-2H)$ where $H$ denotes the mean curvature of the surface $M^2$. The Dirac operator $D$ of $M^2$ defined by the formula

\[ D(\varphi)= e_1 \cdot \nabla^{M^2}_{e_1} \varphi + e_2 \cdot  \nabla^{M^2}_{e_2} \varphi \]

can now be  expressed by the covariant derivative $\nabla^{Y^3}$ and the mean curvature vector:

\[ e_1 \cdot \nabla^{Y^3}_{e_1} (\Phi) + e_2 \cdot \nabla^{Y^3}_{e_2} (\Phi) = D(\varphi) + H \cdot \vec{N} \cdot \varphi . \]

Suppose that the spinor field $\Phi$ on $Y^3$ is a real Killing spinor, i.e. there exists a number $\lambda \in {\Bbb R}^1$ such that for any tangent vector $\vec{T} \in T(Y^3)$ the derivative of $\Phi$ in the direction of $\vec{T}$ is given by the Clifford multiplication:

\[ \nabla^{Y^3}_{\vec{T}} (\Phi) = \lambda \cdot \vec{T} \cdot \Phi . \]

For the restriction $\varphi = \Phi_{|M^2}$ we obtain immediately the equation

\[ D(\varphi) = -2 \lambda \varphi - H \cdot \vec{N} \cdot \varphi . \]

Using the decomposition  $\varphi = \varphi^+ + \varphi^-$ the last equation is equivalent to the  pair of  equations

\[ D(\varphi^+)= (-2 \lambda - iH) \varphi^- \quad , \quad D(\varphi^-)=(-2 \lambda + i H) \varphi^+ . \]

We discuss two special cases.\\

{\bf Proposition 1:} {\it Let $M^2$ be a minimal surface in $Y^3$. Then the restriction  $\varphi = \Phi_{|M^2}$ of any real Killing spinor $\Phi$ on $Y^3$ is an eigenspinor of constant length on the surface $M^2$: 

\[ D(\varphi)= - 2 \lambda \varphi . \] }

On the other hand, suppose that $\Phi$ is a parallel spinor $(\lambda =0)$ on $Y^3$. Then we obtain

\[ D(\varphi^+)= -i H \varphi^- \quad , \quad D(\varphi^-)= iH \varphi^+ . \]

If we introduce the spinor field $\varphi^* = \varphi^+ - i \varphi^-$, a simple calculation shows

\[ D(\varphi^*) =H \varphi^* . \]

The spinor field $\varphi^*$ is given by

\[ \varphi^* = \varphi^+ - i \varphi^- = \frac{1}{2} (\varphi + i \cdot \vec{N} \cdot \varphi) - \frac{i}{2} (\varphi - i \cdot \vec{N} \cdot \varphi) = \frac{1}{2} (1-i) \varphi + \frac{1}{2} (-1+i) \cdot \vec{N} \cdot \varphi . \]

Moreover, the length of $\varphi^*$ is constant. This construction yields the\\

{\bf Proposition 2:} \, {\it Let $\Phi$ be a parallel spinor field defined on the 3-manifold $Y^3$ and denote by $\varphi = \Phi_{|M^2}$ its restriction to $M^2$. Define the spinor field $\varphi^*$ on $M^2$ by the formula

\[ \varphi^* = \frac{1}{2} (1-i) \varphi + \frac{1}{2} (-1+i) \cdot \vec{N} \cdot \varphi . \]

Then $\varphi^*$ is a spinor field of constant length on $M^2$ satisfying the Dirac equation

\[ D(\varphi^*) = H \varphi^* \]

where $H$ denotes  the mean curvature.}\\

{\bf Remark 1:} \, { The map $\Phi \longmapsto \varphi^*$ associating to any parallel spinor $\Phi$ on $Y^3$ a solution of the equation $D(\varphi^*) =H \varphi^*$ is injective.}\\

{\bf Remark 2:} \, { We can apply the above mentioned formulas not only for Killing spinors. Indeed, for any spinor field $\Phi$ we have

\[ D_{Y^3} (\Phi) = D(\varphi) + H \cdot \vec{N} \cdot \varphi + \vec{N} \cdot (\nabla^{Y^3}_{\vec{N}} \Phi) \]

where $D_{Y^3}$ is the Dirac operator of the 3-manifold $Y^3$. Suppose there exists a function $\kappa : M^2 \to {\Bbb C}$ such that the normal derivative $(\nabla^{Y^3}_{\vec{N}} \Phi)$ of the spinor field $\Phi$ is described by $\kappa$:

\[ \left( \nabla^{Y^3}_{\vec{N}} \Phi \right) = \kappa \Phi . \]

Then we obtain

\[ D_{Y^3} (\Phi)= D(\varphi) + (H + \kappa ) \vec{N} \cdot \varphi . \]

This formula (in arbitrary dimension) has been used for the calculation of the spectrum of the Dirac operator on hypersurfaces of the Euclidean space (see [B\"a], [Tr1], [Tr2]).\\

\section{Solutions of the Dirac Equation with Potential on Riemannian Surfaces.}

Let $(M^2,g)$ be an oriented, 2-dimensional Riemannian  manifold with spin structure. $H:M^2 \to {\Bbb R}^1$ denotes a given smooth, real-valued function defined on the surface. In this part we study spinor fields $\varphi$ on $M^2$ that are solutions of the differential equation

\[ D(\varphi)= H \varphi . \]

If  we decompose the spinor field into $\varphi = \varphi^+ + \varphi^-$ according to the splitting $S=S^+ \oplus S^-$ of the spinor bundle  the equation we want to study is equivalent to the system 

\[ D(\varphi^+)= H \varphi^- \quad , \quad D(\varphi^-)= H \varphi^+ . \]

To any solution $\varphi$ of this equation we associate two forms $F_{\pm}$ defined for pairs $X,Y \in T(M^2)$ of tangent vectors: 

\[ F_+ (X,Y)= \mbox{Re} (\nabla_X \varphi^+ , Y \cdot \varphi^-) \quad , \quad  F_- (X,Y)= \mbox{Re} (\nabla_X \varphi^- , Y \cdot \varphi^+) . \]

{\bf Proposition 3:} \, {\it
\begin{itemize}
\item[a.)] $F_{\pm}$ are symmetric bilinear forms on $T(M^2)$.
\item[b.)] The trace of $F_{\pm}$ is given by 
$Tr(F_{\pm})=-H| \varphi^{\mp}|^2 $.
\end{itemize}}

{\bf Proof:} The symmetry of $F_{\pm}$ is a consequence of the Dirac equation as well as the assumption that $H$ is a real-valued function. Indeed, we have 

\begin{eqnarray*}
\mbox{Re} (\nabla_{e_1} \varphi^+, e_2 \varphi^-)&=& \mbox{Re} (e_1 \cdot \nabla_{e_1} \varphi^+, e_1 \cdot e_2 \cdot \varphi^-)=  \mbox{Re} (H \varphi^- - e_2 \cdot \nabla_{e_2} \varphi^+, e_1 \cdot e_2 \cdot \varphi^-) =\\
&=& H \cdot \mbox{Re} (\varphi^- , e_1 \cdot e_2 \cdot \varphi^-) + \mbox{Re} (\nabla_{e_2} \varphi^+ , e_2 \cdot e_1 \cdot e_2  \cdot \varphi^-)=\\
&=&  0 + \mbox{Re} (\nabla_{e_2} \varphi^+ , e_1 \cdot \varphi^-) . 
\end{eqnarray*}

Moreover, we calculate the trace of $F_{\pm}$:

\[
Tr (F_{\pm})= \mbox{Re} (\nabla_{e_1} \varphi^{\pm}, e_1 \cdot \varphi^{\mp}) + \mbox{Re} (\nabla_{e_2} \varphi^{\pm}, e_2 \cdot \varphi^{\mp})  = - \mbox{Re} (D(\varphi^{\pm}), \varphi^{\mp})= - H|\varphi^{\mp}|^2 . 
\]
\mbox{} \hfill \rule{3mm}{3mm}\\

We study now special solutions of the equation $D(\varphi)=H \varphi$, i.e. solutions with {\it constant length} $|\varphi | \equiv const \not= 0$. It may happen that the components $\varphi^{\pm}$ have a non-empty zero set.\\

{\bf Proposition 4:} \, {\it Suppose that the spinor field $\varphi$ defined on the Riemannian surface $M^2$ is a solution of the equation

\[ D(\varphi)=H \varphi \hspace{1cm} \mbox{with} \,\,\,\, |\varphi | \equiv const \not= 0 . \]

Then the forms $F_{\pm}$ are related by the equation}

\[ | \varphi^+ |^2 F_+ = | \varphi^- |^2 F_- . \]

{\bf Proof:} In case one of the spinors $\varphi^+$ or $\varphi^-$ vanishes at a fixed point $m_0 \in M^2$ the relation between $F_+$ and $F_-$ is trivial. Otherwise there exists a neighbourhood $V$ of the point $m_0 \in M^2$ such that both spinors $\varphi^+$ and $\varphi^-$ are not zero at any point $m \in V$. The spinors

\[ \frac{e_1 \cdot \varphi^-}{|\varphi^-|} \quad , \quad \frac{e_2 \cdot \varphi^-}{|\varphi^-|} \]

are an orthonormal base in $S^+$ with respect to the Euclidean scalar product $\mbox{Re} ( \cdot , \cdot )$. Therefore we obtain (on $V$)

\begin{eqnarray*}
\nabla_X \varphi^+ &=& \mbox{Re} \left( \nabla_X \varphi^+, \frac{e_1 \cdot 
\varphi^-}{|\varphi^- |} \right) \frac{e_1 \cdot \varphi^-}{|\varphi^- |} 
 + \mbox{Re} \left( \nabla_X \varphi^+ , \frac{e_2 \cdot \varphi^-}{|\varphi^-|} \right) \frac{e_2 \cdot \varphi^-}{|\varphi^- | } =\\
&=& \frac{1}{| \varphi^- |^2} \{ F_+ (X, e_1) e_1 + F_+ (X, e_2) e_2 \} \cdot \varphi^- . 
\end{eqnarray*}

A similar calculation yields the formula 

\[ \nabla_X \varphi^- = \frac{1}{| \varphi^+ |^2} \{ F_- (X, e_1) e_1 + F_- (X, e_2) e_2 \} \cdot \varphi^+ . \]

We multiply the equations by $\varphi^+$ and $\varphi^-$ respectively and sum up. Then we obtain

\[ \frac{1}{2} \nabla_X (| \varphi^+ |^2 + | \varphi^- |^2 )= \mbox{Re} (A (X) \varphi^- , \varphi^+ ) \]

where the endomorphism $A:T(M^2) \to T(M^2)$ is defined by

\[
A(X) = \left\{ \frac{F_+ (X, e_1)}{| \varphi^- |^2} - \frac{F_- (X, e_1)}{| \varphi^+ |^2} \right\} e_1 + \left\{ \frac{F_+ (X, e_2 )}{| \varphi^- |^2} - \frac{F_- (X, e_2 )}{ |\varphi^+ |^2 } \right\} e_2 . 
\]

Since $F_{\pm}$ are symmetric tensors, the endomorphism $A$ is symmetric too. Moreover, the trace of  $A$ vanishes:

\[ Tr A= \frac{1}{| \varphi^- |^2} Tr(F_+ ) - \frac{1}{|€\varphi^+ |^2} Tr (F_-)= - H + H =0 . \]

The length of the spinor field $\varphi$ is constant. This implies

\[ \mbox{Re} (A(X) \cdot \varphi^- , \varphi^+ )=0 . \]

At any point $m \in V$ of the set $V$ the spinors $\varphi^+ , \varphi^- $ are non-trivial. Then  the rank of the endomorphisms $A:T(M^2) \to T(M^2)$ is not greater than one. All in all, $A$ is symmetric, $Tr (A)=0$ and $rg (A) \le 1$, i.e. $A \equiv 0$. \\
\mbox{} \hfill \rule{3mm}{3mm}\\

We now consider the sum

\[ F=F_+ + F_- . \]

At points with $\varphi^+ \not= 0$ (or $\varphi^- \not= 0$) we have

\[ \frac{F}{|\varphi|^2} = \frac{F_+ + F_- }{|\varphi^+ |^2 + |\varphi^- |^2} =  \frac{\left( \frac{|\varphi^-|^2}{|\varphi^+|^2} + 1 \right) F_-}{|\varphi^+ |^2 + |\varphi^- |^2} = \frac{F_-}{|\varphi^+ |^2} \]

as well as 

\[ \frac{F}{|\varphi|^2} = \frac{F_+ + F_- }{|\varphi^+ |^2 + |\varphi^- |^2} =  \frac{\left(   1 + \frac{|\varphi^+|^2}{|\varphi^-|^2}\right) F_+}{|\varphi^+ |^2 + |\varphi^- |^2} = \frac{F_+}{|\varphi^- |^2} .  \]

The endomorphism $E:T(M^2) \to T(M^2)$, $E= \frac{F}{|\varphi|^2}$, is defined at all points of $M^2$ and the formulas derived in the proof of proposition 4 in fact prove the following \\

{\bf Proposition 5:} \, {\it Let $\varphi$ be a solution of the differential equation $D(\varphi) = H \varphi $ on a Riemannian surface $(M^2,g)$ with a real-valued function $H:M^2 \to {\Bbb R}^1$. Suppose that the length $|\varphi | \equiv const \not= 0$ of the spinor field $\varphi$ is constant. Then}

\[ E(X)= \frac{1}{|\varphi|^2} \mbox{Re} (\nabla_X \varphi , Y \cdot \varphi ) \]

{\it is a symmetric endomorphism $E:T(M^2) \to T(M^2)$ such that}
\begin{itemize}
\item[a.)] $\nabla_X \varphi^+ =E (X) \cdot \varphi^- \quad , \quad  \nabla_X \varphi^- =E (X) \cdot \varphi^+$
\item[b.)] $Tr (E)= - H$.
\end{itemize}
\mbox{} \hfill \rule{3mm}{3mm}\\

For a given triple $(M^2,g,E)$ of a Riemannian surface and symmetric endomorphism the existence of a non-trivial solution $\varphi$ of the equation

\[ \nabla_X \varphi = E(X) \cdot \varphi \]

implies certain integrability conditions. It turns out that in this way we obtain  precisely the well-known Gau{\ss} and Codazzi equations of the classical theory of surfaces in Euclidean 3-space.\\

{\bf Proposition 6:} \, {\it Let $(M^2,g)$ be a 2-dimensional  Riemannian surface with a fixed spin structure and suppose that $E:T(M^2) \to T(M^2)$ is a symmetric endomorphism. If there exists a non-trivial solution of the equation}

\[ \nabla_X \varphi =E (X) \cdot \varphi \quad , \quad X \in T(M^2) \]

{\it then
\begin{itemize}
\item[a.)] (Codazzi equation): 
\quad $\nabla_X (E(Y)) - \nabla_Y (E(X)) - E([X,Y]) =0 .$
\item[b.)] (Gau{\ss} equation): 
\quad $\det (E) = \frac{1}{4} G$, \, \,  where $G$ is the Gaussian curvature 
of $(M^2,g)$.
\end{itemize}}

{\bf Proof:} We prove the two equations in a way similar to the derivation of  the integrability conditions for the Riemannian metric in case the space admits a Killing spinor (see [BFGK]). We differentiate the equation

\[ \nabla_X \varphi =E(X) \cdot \varphi \]

and then we calculate the curvature tensor $R^S$ of the spinor bundle $S$:

\begin{eqnarray*}
R^S (X,Y) \varphi &=& \nabla_X \nabla_Y \varphi - \nabla_Y \nabla_X \varphi - \nabla_{[X,Y]} \varphi=\\
&=& \{ \nabla_X (E(Y)) - \nabla_Y (E(X)) - E([X,Y])  + E(Y) E(X) - E(X)E(Y) \} \cdot \varphi .
\end{eqnarray*}

On the other side, the curvature tensor $R^S :S \to S$ is given by the formula

\[ R^S (e_1, e_2 )= \frac{1}{2} R_{1212} \, \, e_1 \cdot e_2 . \]

Denote by $A(X,Y)$ the differential of $E$:

\[ A(X,Y) = \nabla_X (E(Y)) - \nabla_Y (E(X)) - E([X,Y]) . \]

A simple algebraic calculation in the spin representation then leads to the equations

\[ - A(X,Y) \varphi^- = \left( 2 \mbox{det} (E) + \frac{R_{1212}}{2} \right) i \varphi^+ \]

\[ A(X,Y) \varphi^+ = \left( 2 \mbox{det} (E) + \frac{R_{1212}}{2} \right) i \varphi^- . \]

We multiply the first equation once by the vector $A(X;Y)$:

\[ ||A(X,Y)||^2 \varphi^- = - \left(2 \mbox{det} (E) + \frac{R_{1212}}{2} \right)^2 \varphi^- \]

and then we conclude $A(X,Y) \equiv 0$ (Codazzi equation) as well as  $\mbox{det} (E) = - \frac{1}{4} R_{1212} = \frac{1}{4} G$ (Gau{\ss} equation). \\
\mbox{} \hfill \rule{3mm}{3mm}\\

For a given triple $(M^2,g,E)$ consisting of a Riemannian spin surface $(M^2,g)$ and of a symmetric endomorphism $E$ we will denote by ${\cal K} (M^2,g,E)$ the space of all spinor fields $\varphi$ satisfying the equation $\nabla_X \varphi =E(X) \cdot \varphi$. It is invariant under the quaternionic structure $\alpha :S \to S$, i.e. ${\cal K} (M^2,g,E)$ is a quaternionic vector space (see $\S 4$). Denote by $(-H)$ the trace of $E$, 

\[ Tr (E)= - H . \]

Then we have

\[ {\cal K} (M^2,g,E) \subset \ker (D - H) . \]

In this part of the paper we prove that any spinor field $\varphi \in \ker (D - H)$ of constant length belongs to one of the subspaces ${\cal K} (M^2,g,E)$ for a suitable symmetric endomorphism $E$, $Tr (E) = - H$. \\


Finally, we consider the lengths

\[ L_+ = ||€\varphi^+ ||^2 \quad , \quad L_- =||€\varphi^- ||^2 \]

of a non-trivial solution $\varphi \in {\cal K} (M^2,g,E)$.  Using the integrability condition $\mbox{det} \, (E) = \frac{G}{4}$ (i.e. $||E||^2 =H^2 - \frac{G}{2}$) as well as the well-known formula $D^2 =\Delta + \frac{G}{2}$ for the square $D^2$ of the Dirac operator we can derive formulas for $\Delta (L_{\pm})$:

\begin{eqnarray*}
\Delta (L_{\pm}) &=& 2 (\Delta (\varphi^{\pm}), \varphi^{\pm}) - 2 \langle \nabla (\varphi^{\pm}), \nabla (\varphi^{\pm}) \rangle\\
&=& 2 (D^2 (\varphi^{\pm}), \varphi^{\pm} ) - 2 \left( \frac{G}{2} \right) \bullet ||\varphi^{\pm}||^2  -  2 ||E||^2 || \varphi^{\mp} ||^2\\
&=& 2 \left( H^2 - \frac{G}{2} \right) (L_{\pm} - L_{\mp}) + 2 \mbox{Re} (grad \, (H) \cdot \varphi^{\mp}, \varphi^{\pm})
\end{eqnarray*}

In particular, if $H \equiv const$ is constant,  the difference $u= L_+ - L_-$ satisfies the differential equation

\[ \Delta (u) = 4 \left( H^2 - \frac{G}{2} \right) u . \]

\section{The Period Form of a Spinor with $\nabla_X \varphi =E(X) \cdot \varphi$.}

We consider  a spinor field $\varphi$ on a Riemannian surface $(M^2,g)$ such that

\[ \nabla_X \varphi =E(X) \cdot \varphi \]

for a fixed symmetric endomorphism $E$. The spinor bundle $S$ carries a quaternionic structure $\alpha : S \to S$ commuting with Clifford multiplication and interchanging the decomposition $S= S^+ \oplus S^-$ (see [F]). For any spinor field $\varphi = \varphi^+ + \varphi^-$ we define three 1-forms by

\[ \xi^{\varphi} (X) = 2 (X \cdot \varphi^+, \varphi^- ) \]

\[ \xi^{\varphi}_+ (X) =(X \cdot \varphi^+, \alpha (\varphi^+)) \quad , \quad \xi^{\varphi}_- (X)=(X \cdot \varphi^-, \alpha (\varphi^-)) . \]

$\xi^{\varphi}$ and $\xi^{\varphi}_+$ are $\Lambda^{1,0}$-forms,  $\xi^{\varphi}_-$ is a $\Lambda^{0,1}$-form. Indeed, $e_1 \cdot e_2$ acts on $S^+$ (on $S^-$) by multiplication by $(-i)$ (by $i$). Now we obtain

\[
(\star \xi^{\varphi} )(e_1) = - \xi^{\varphi} (e_2) =2 (-e_2 \cdot \varphi^+, \varphi^-)=(-i e_2 \cdot e_1 \cdot e_2 \cdot \varphi^+, \varphi^-)= -i \xi^{\varphi} (e_1) , 
\]

i.e. $\star \xi^{\varphi} =-i \xi^{\varphi}$ holds. A similar calculation gives $\star \xi^{\varphi}_+ =-i \xi^{\varphi}_+$ and $\star \xi^{\varphi}_- = i \xi^{\varphi}_-$. We split the 1-form $\xi^{\varphi}$ into its real and imaginary part:

\[ \xi^{\varphi} = w^{\varphi} + i \mu^{\varphi} . \]

Moreover, we introduce the 1-form $\Omega^{\varphi}$

\[ \Omega^{\varphi} = \xi^{\varphi}_+ - \xi^{\varphi}_- . \]

Then we have\\

{\bf Proposition 7:} \, {\it Let $(M^2,g)$ be a Riemannian spin surface and $E:T(M^2) \to  T(M^2)$ a symmetric endomorphism of trace $-H$. Suppose the spinor field $\varphi$ is a solution of the equation $ \nabla_X \varphi =E(X) \cdot \varphi $. Then}
\begin{itemize}
\item[a.)] $dw^{\varphi} =0.$
\item[b.)] $d \mu^{\varphi} =2H \{ |\varphi^- |^2 - |\varphi^+ |^2 \} dM^2.$
\item[c.)] $d \Omega^{\varphi} =0.$
\end{itemize}

{\bf Proof:} We calculate $dw^{\varphi}$:

\begin{eqnarray*}
\frac{1}{2} dw^{\varphi} (X,Y) &=& X(\mbox{Re}(Y \cdot \varphi^+, \varphi^- )) - Y(\mbox{Re} (X \cdot \varphi^+, \varphi^-)) - \mbox{Re} ([X,Y] \cdot \varphi^+, \varphi^-)=\\
&=& \{ g(X,E(Y)) - g(Y,E(X)) \} |\varphi^- |^2   + \{ g(X,E(Y)) - g(Y,E(X)) \} |\varphi^+ |^2  .
\end{eqnarray*}

Since $E$ is symmetric, we obtain $dw^{\varphi}=0$. A similar calculation shows the formula for $d \mu^{\varphi}$. For the proof of  $d \Omega^{\varphi}=0$ we first remark that the quaternionic structure $\alpha :S \to S$ and the hermitian product $( \cdot , \cdot)$ on $S$ are related by

\[ (\varphi_1, \alpha (\varphi_2)) = - (\overline{\alpha (\varphi_1), \varphi_2}) . \]

Using this formula we can transform $d \xi^{\varphi}_-$ in the following way:

\begin{eqnarray*}
d \xi^{\varphi}_- (X,Y) &=& (Y \cdot E(X) \cdot \varphi^+, \alpha (\varphi^-)) + (Y \cdot \varphi^-, \alpha (E(X) \cdot \varphi^+))\\
& & - (X \cdot E(Y) \cdot \varphi^+, \alpha (\varphi^-)) - (X \cdot \varphi^-, \alpha (E(Y) \cdot \varphi^+))\\
&=& - (\overline{\alpha (Y \cdot E(X) \cdot \varphi^+), \varphi^-}) - (E(X) \cdot Y \cdot \varphi^-, \alpha (\varphi^+))\\
& &+ ( \overline{\alpha (X \cdot E(Y) \cdot \varphi^+), \varphi^-}) + (E(Y) \cdot X \cdot \varphi^-, \alpha (\varphi^+))\\
&=& - (E(X) \cdot Y \cdot \varphi^-, \alpha (\varphi^+)) - (E(X) \cdot Y \cdot \varphi^-, \alpha (\varphi^+))\\
& &+ (E(Y) \cdot X \cdot \varphi^-, \alpha (\varphi^+)) + (E(Y) \cdot X \cdot \varphi^-, \alpha (\varphi^+)) . 
\end{eqnarray*}

On the other hand we calculate $d \xi^{\varphi}_+$:

\begin{eqnarray*}
d \xi^{\varphi}_+ (X,Y) &=& (Y \cdot E(X) \cdot \varphi^-, \alpha (\varphi^+)) + (Y \cdot \varphi^+, \alpha (E(X) \cdot \varphi^-))\\
& & - (X \cdot E(Y) \cdot \varphi^-, \alpha (\varphi^+)) - (X \cdot \varphi^+, \alpha (E(Y) \cdot \varphi^-))=\\
&=& (Y \cdot E(X) \cdot \varphi^-, \alpha (\varphi^+)) - (E(X) \cdot Y \cdot \varphi^+, \alpha (\varphi^-)) \\
& & - (X \cdot E(Y) \cdot \varphi^-, \alpha (\varphi^+)) + (E(Y) \cdot X \cdot \varphi^+, \alpha (\varphi^-)). 
\end{eqnarray*}

Finally we obtain

\begin{eqnarray*}
d(\xi^{\varphi}_- - \xi^{\varphi}_+)(X,Y) &=& - (\{ E(X) \cdot Y + Y \cdot E(X) \} \varphi^-, \alpha (\varphi^+))\\
& & + (\{ E(Y) \cdot X + X \cdot E(Y)\} \varphi^-, \alpha (\varphi^+)) - \\
&=& 2 \{ g (E(X), Y) - g (E(Y), X) \} (\varphi^-, \alpha (\varphi^+ ))
\end{eqnarray*}

and $d(\xi^{\varphi}_- - \xi^{\varphi}_+)=0$ follows again by the symmetry of $E$.\\
\mbox{} \hfill \rule{3mm}{3mm}\\

Let us consider the case that $(M^2,g)$ is isometrically immersed into the Euclidean space ${\Bbb R}^3$, $\Phi$ is a parallel spinor on ${\Bbb R}^3$ and the spinor field $\varphi^*$ on $M^2$ defined by the formula

\[ \varphi^* = \frac{1}{2} \left( \Phi_{|M^2} + i \cdot \vec{N} \cdot \Phi_{|M^2} \right) + \frac{i}{2} \left( i \cdot \vec{N} \cdot \Phi_{|M^2} - \Phi_{|M^2} \right) \]

(see $\S 2$). In this case the forms $w^{\varphi^*}$ and $\Omega^{\varphi^*}$ are given by the expressions 

\[ w^{\varphi^*} (X) = - \mbox{Im} (X \cdot \Phi, \Phi) \quad , \quad \Omega^{\varphi^*} (X) =(X \cdot \Phi, \alpha (\Phi)) , \]

and  are exact 1-forms. Indeed, we defined functions $f: {\Bbb R}^3 \to {\Bbb R}^1$ and $g: {\Bbb R}^3 \to {\Bbb C}$ by \[ f(m) = - \mbox{Im} \langle m \cdot \Phi, \Phi \rangle \quad , \quad g (m)= \langle m \cdot \Phi, \alpha (\Phi) \rangle \]

and then we have $df= w^{\varphi^*}, dg= \Omega^{\varphi^*}$. We remark that $f$ and $g$ describe in fact the isometric immersion $M^2 \hookrightarrow {\Bbb R}^3$ we started with. The 3-dimensional spinor $\Phi \in \Delta_3$ defines a real 3-dimensional subspace $\Delta_3 (\Phi)$ by

\[ \Delta_3 (\Phi)= \{ \Psi \in \Delta_3: \quad \mbox{Re} (\Psi, \Phi)=0 \} . \]

The map $\Psi \to (- \mbox{Im} (\Psi, \Phi), (\Psi, \alpha (\Phi)))$ is an isometry between $\Delta_3 (\Phi)$ and ${\Bbb R}^1 \oplus {\Bbb C} = {\Bbb R}^3$. Clearly, the immersion $M^2 \hookrightarrow {\Bbb R}^3$ is given by

\[ M^2 \ni m \longmapsto m \cdot \Phi \in \Delta_3 (\Phi) , \]

i.e. by the functions $f_{|M^2}$ and $g_{|M^2}$. With respect to $d(f_{|M^2})=w^{\varphi^*}$ and $d(g_{|M^2})= \Omega^{\varphi^*}$ we obtain a formula for the isometric immersion $M^2 \hookrightarrow {\Bbb R}^3$:

\[ \oint (w^{\varphi^*}, \Omega^{\varphi^*}) : M^2 \to {\Bbb R}^3 . \]

(Weierstra{\ss} representation of the surface.)\\

In general, we call a solution $\varphi$ of the differential equation $\nabla_X \varphi =E(X) \cdot \varphi$ {\it exact} iff the corresponding forms $w^{\varphi}, \Omega^{\varphi}$ are exact 1-forms. Using the definition

\[ Hess \,(h)(X,Y) = \frac{1}{2} \{ g (\nabla_X (grad \, (h)), Y)  + g (X, \nabla_Y (grad \, (h))) \} \]

of the Hessian of a smooth function $h$ defined on a Riemannian manifold we obtain the \\

{\bf Proposition 8:} \, {\it  Let $\varphi \in {\cal K} (M^2, g,E)$ be an exact solution of the differential equations $\nabla_X \varphi =E (X) \cdot  \varphi$ with $df= w^{\varphi}$, $dg = \Omega^{\varphi}$. Then}

\begin{itemize}
\item[a.)] $Hess \, (f)=2 \left( |\varphi^+ |^2 - | \varphi^- |^2 \right) E.$
\item[b.)] $|grad \, f|^2 = 4 | \varphi^+ |^2 | \varphi^- |^2$.
\item[c.)] $Hess \, (g) = - 4 ( \varphi^-, \alpha (\varphi^+ )) E$.
\item[d.)] $|grad \, (g)|^2 = ( |\varphi^+ |^2 - |\varphi^- |^2 )^2$. 
\end{itemize}
\mbox{} \hfill \rule{3mm}{3mm}\\

In particular, the determinant of the Hessian of the function $f$ is given by

\[ \mbox{det} (Hess \, (f))=4 (|\varphi^+ |^2 - |\varphi^- |^2 )^2 \mbox{det} (E)=(|\varphi^+ |^2 - |\varphi^- |^2 )^2 G . \]

Here we used proposition 6, i.e. $\mbox{det} (E)= \frac{1}{4} G$.\\

{\bf Corollary:} \, {\it Let $M^2$ be a compact Riemannian spin-manifold and suppose that $\varphi \in {\cal K} (M^2,g,E)$ is an exact, non-trivial solution. Then the spinors $\varphi^+$ or $\varphi^-$ vanish at least at  out one point.  Moreover, there exists $m_0 \in M^2$ such that $G(m_0) \ge 0$.}\\

{\bf Proof:} At a maximum point $m_0 \in M^2$ of $f$ we have

\[ grad (f) (m_0) =0 \quad , \quad \mbox{det} (Hess \, (f) (m_0)) \ge 0 . \]
\mbox{} \hfill \rule{3mm}{3mm}\\

Recall that for any 2-dimensional Riemannian manifold $(M^2,g)$ and any function $h :M^2 \to {\Bbb R}^1$ the 2-form

\[ \{ 2 \det (Hess \, (h)) - |grad (h)|^2 G \} dM^2 = d \mu^1 \]

is exact (see [S], page 47). Using this formula in case of an exact solution $\varphi \in {\cal K} (M^2,g,E)$ we obtain

\[ \int\limits_{M^2} (|\varphi^+ |^2 - |\varphi^- |^2 )^2 G= 2 \int\limits_{M^2} |\varphi^+ |^2 |\varphi^- |^2 G . \]

{\bf Corollary:} \, {\it Let $M^2$ be a compact Riemannian spin manifold and suppose that $\varphi \in {\cal K} (M^2,g,E)$ is an exact solution. Then }

\[  \int\limits_{M^2} (|\varphi|^4 - 6 |\varphi^+|^2 |\varphi^-|^2 )G=0 . \]
\mbox{} \hfill \rule{3mm}{3mm}\\

We again discuss the last formula in case of an isometrically immersed surface $M^2 \hookrightarrow {\Bbb R}^3$ and a given parallel spinor $\Phi$ on ${\Bbb R}^3$. We apply the integral  formula to the spinor $\varphi^* = \varphi^*_+ + \varphi^*_-$ where

\[ \varphi^*_+ = \frac{1}{2} \left( \Phi + i \vec{N} \cdot \Phi \right) \quad , \quad \varphi^*_- = \frac{1}{2} \left( - \Phi + i \cdot \vec{N} \cdot \Phi \right) . \]

In this case we have

\[ |\varphi^*_+ |^2 = \frac{1}{2} |\Phi|^2 + \frac{1}{2} \langle i \vec{N} \cdot \Phi, \Phi \rangle \quad , \quad  |\varphi^*_- |^2 = \frac{1}{2} |\Phi|^2 - \frac{1}{2} \langle i \vec{N} \cdot \Phi, \Phi \rangle \]

and $(|\Phi|\equiv 1)$ therefore

\[ 1 - 6 |\varphi^*_+ |^2 |\varphi^*_-|^2 = - \frac{1}{2} + \frac{3}{2} \langle i \cdot \vec{N} \cdot \Phi, \Phi \rangle^2 . \]

Consequently the integral formula yields

\[ \int\limits_{M^2} G= 3 \int\limits_{M^2} \langle i \vec{N} \Phi, \Phi \rangle^2 G . \]

\newfont{\graf}{eufm10}
\newcommand{\alta}{\mbox{\graf a}}

The spinors $i \Phi$ as well as $\vec{N} \cdot \Phi$ belong to $V(\Phi) \subset \Delta_3$, the space of the immersion $M^2 \hookrightarrow {\Bbb R}^3 =V ( \Phi)$. The last formula means therefore

\[ \int\limits_{M^2} G= 3 \int\limits_{M^2} \langle \vec{N}, {\alta}_3 \rangle^2 G \]

for the unit vector ${\alta}_3 = i \Phi \in V(\Phi) = {\Bbb R}^3$. 

\section{The Spin Formulation of the Theory of Surfaces in ${\Bbb R}^3$.}

An oriented, immersed surface $M^2 \hookrightarrow {\Bbb R}^3$ inherits from ${\Bbb R}^3$ an inner metric $g$, a spin structure and a  solution $\varphi$ of the Dirac equation

\[ D(\varphi)= H \varphi \]

of constant length $|\varphi | \equiv 1$ where $H$ denotes the mean curvature of the surface. The spinor field $\varphi$ on $M^2$ is the restriction of a parallel spinor field $\Phi$ of the Euclidean space ${\Bbb R}^3$. The period forms $w^{\varphi}$ and $\Omega^{\varphi}$ are exact and the immersion $M^2 \hookrightarrow {\Bbb R}^3$ is given by integration of the ${\Bbb R}^1 \oplus {\Bbb C} = {\Bbb R}^3$ valued form $(w^{\varphi}, \Omega^{\varphi})$. At least  locally the converse is true: Given an oriented, 2-dimensional Riemannian manifold $(M^2,g)$ with a fixed spin structure and a solution of constant length of the Dirac equation $D(\varphi )= H \varphi$ for some smooth function $H:M^2 \to {\Bbb R}^1$, there exists a symmetric endomorphism $ E:T(M^2) \to T(M^2)$ such that $\varphi \in {\cal K} (M^2,g,E)$. Moreover, $2 E$ is the second fundamental form of an isometric immersion $(M^2,g) \to {\Bbb R}^3$. We formulate this description of the theory of surfaces in ${\Bbb R}^3$ in the following\\

{\bf Theorem 1:} \, {\it Let $(M^2,g)$ be an oriented, 2-dimensional Riemannian manifold and $H:M^2 \to {\Bbb R}^1$ a smooth function. Then there is a correspondence between the following datas:
\begin{enumerate}
\item an isometric immersion $(\tilde{M}^2, g) \to {\Bbb R}^3$ of the universal covering $\tilde{M}^2$ into the Euclidean space ${\Bbb R}^3$ with mean curvature $H$. 
\item a solution $\varphi$ with constant length $|\varphi | \equiv 1$ of the Dirac equation $D(\varphi)=H \cdot \varphi$.
\item a pair $(\varphi, E)$ consisting  of a symmetric endomorphism $E$ such 
that \linebreak $Tr(E)=-H$ and a spinor field $\varphi$ satisfying the 
equation $\nabla_X \varphi =E(X) \cdot \varphi$.
\end{enumerate}}
\mbox{} \hfill \rule{3mm}{3mm}\\

We apply now the well-known formulas for the change of the Dirac operator under a conformal change of the metric. Suppose that $\tilde{g} = \sigma g$ are two conformally equivalent metrics on $M^2$ where $\sigma :M^2 \to (0, \infty)$ is a positive function. Denote by $D$ and $\tilde{D}$ the Dirac operator corresponding to the metric $g$ and $\tilde{g}$ respectively. Then

\[ \tilde{D} (\tilde{\varphi}) = \sigma^{-3/4} \widetilde{D(\sigma^{1/4}
 \varphi)} \]

holds (see [BFGK]). Let us consider a solution $\varphi$ of the Dirac equation

\[ D(\varphi)= \lambda \varphi \]

on $(M^2,g)$ and suppose that $\varphi$ never vanishes. We introduce the Riemannian metric $\tilde{g} = |\varphi |^4 g$ as well as the spinor field 
$\varphi^* = \frac{\varphi}{|\varphi |}$. Then we obtain

\[ \tilde{D} (\varphi^*) = \frac{\lambda}{|\varphi |^2} \varphi^* \quad , \quad |\varphi^* | \equiv 1 , \]

and thus an isometric immersion $(\tilde{M}^2, |\varphi |^4 g) \to {\Bbb R}^3$ with mean curvature $H= \frac{\lambda}{|\varphi |^2}$.\\

{\bf Theorem 2:} \, {\it Let $(M^2,g)$ be an oriented, 2-dimensional 
Riemannian manifold. Any spinor field $\varphi$ without zeros that is a 
solution of the equation}

\[ D(\varphi) = \lambda \varphi \]

{\it defines an isometric immersion $(\tilde{M}^2, |\varphi |^4 g) 
\hookrightarrow {\Bbb R}^3$ with mean curvature 
$H= \frac{\lambda}{|\varphi |^2}$. }\\
\mbox{} \hfill \rule{3mm}{3mm}\\

{\bf Remark 1:} Consider the case that $M^2 \hookrightarrow S^3$ is a minimal surface in $S^3$. Let $\Phi$ be a real Killing spinor on $S^3$, i.e.

\[ \nabla_{\vec{T}} (\Phi)= \frac{1}{2} \vec{T} \cdot \Phi . \]

The restriction $\varphi = \Phi_{|M^2}$ is an eigenspinor of the Dirac operator on $M^2$ with constant length (proposition 1). Therefore $\varphi$ defines an isometric immersion of  $(\tilde{M}^2,g) \hookrightarrow {\Bbb R}^3$ with mean curvature $H \equiv -1$. This transformation associates to any minimal surface $M^2 \hookrightarrow S^3$ a surface of constant mean curvature $H \equiv -1$ in ${\Bbb R}^3$, a well-known construction (see [La]).\\

{\bf Remark 2:} Using the described correspondence between isometric immersions of surfaces into ${\Bbb R}^3$ and solutions of the Dirac equation $D(\varphi)=H \cdot \varphi$ one can immediately remark that several statements of the elementary theory of surfaces are equivalent to several statements concerning solutions of the twistor equation (see [BFGK]). For example, in [Li] (see also proposition 6) one can find the following theorem: if $f:M^2 \to {\Bbb R}^1$ is a real-valued function such that the equation

\[ \nabla_{\vec{T}} (\varphi) + \frac{1}{2} f \cdot \vec{T} \cdot \varphi =0 \]

admits a non-trivial solution then $f$ is constant and $f^2=G$. In the theory of surfaces this statement correspondends to the fact that an umbilic surface is a part of the sphere or the plane. Indeed, an umbilic surface $M^2 \hookrightarrow {\Bbb R}^3$ admits a spinor field $\varphi$ such that

\[\nabla_{\vec{T}} (\varphi) + \frac{1}{2} H \vec{T} \cdot \varphi =0 \]

and therefore $H^2 =G = const$, i.e. the second fundamental form is proportional to the metric. In a similar way one can translate other facts of the theory of surfaces into properties of solutions of the equation 
$\nabla_X \varphi =E(X) \cdot \varphi$.\\

\vspace{2cm}
Thomas Friedrich\\
Humboldt-Universit\"at zu Berlin\\
Institut f\"ur Reine Mathematik\\
Sitz: Ziegelstra{\ss}e 13a\\
D-10099 Berlin\\
e-mail: friedric@mathematik.hu-berlin.de

\end{document}